\documentclass[11pt]{article}
\usepackage{amsmath}
\usepackage{amsthm}
\usepackage{amsfonts}
\usepackage{amssymb}
\newcommand{\be}{\begin{equation}}
\newcommand{\ee}{\end{equation}}
\newcommand{\bea}{\begin{eqnarray}}
\newcommand{\eea}{\end{eqnarray}}
\newcommand{\bean}{\begin{eqnarray*}}
\newcommand{\eean}{\end{eqnarray*}}
\newcommand{\brray}{\begin{array}}
\newcommand{\erray}{\end{array}}
\newcommand{\ben}{\begin{equation}{nonumber}}
\newcommand{\een}{\end{equation}{nonumber}}

\newtheorem{dfn}{Definition}[section]
\newtheorem{thm}[dfn]{Theorem}
\newtheorem{lema}[dfn]{Lemma}
\newtheorem{ppsn}[dfn]{Proposition}
\newtheorem{coro}[dfn]{Corollary}
\newtheorem{xmpl}[dfn]{Example}
\newtheorem{rmrk}[dfn]{Remark}

\newcommand{\bdfn}{\begin{dfn}}
\newcommand{\bthm}{\begin{thm}}
\newcommand{\blema}{\begin{lema}}
\newcommand{\bppsn}{\begin{ppsn}}
\newcommand{\bcoro}{\begin{coro}}
\newcommand{\bxmpl}{\begin{xmpl}}
\newcommand{\brmrk}{\begin{rmrk}}

\newcommand{\edfn}{\end{dfn}}
\newcommand{\ethm}{\end{thm}}
\newcommand{\elema}{\end{lema}}
\newcommand{\eppsn}{\end{ppsn}}
\newcommand{\ecoro}{\end{coro}}
\newcommand{\exmpl}{\end{xmpl}}
\newcommand{\ermrk}{\end{rmrk}}

\def\wZ{\widetilde{\mathcal Z}}

\newcommand{\A}{\mathcal A}
\newcommand{\E}{\mathcal E}
\newcommand{\C}{\mathcal C}
\newcommand{\G}{\mathcal G}
\newcommand{\I}{\mathcal I}

\newcommand{\M}{\mathcal M}

\newcommand{\Z}{\mathcal Z}
\newcommand{\MBBC}{\mathbb{C}}
\newcommand{\MBBZ}{\mathbb{Z}}
\newcommand{\MBBR}{\mathbb{R}}
\newcommand{\MBFH}{\mathbf{h}}
\newcommand{\MBFK}{\mathbf{k}}
\newcommand{\EXP}{\textbf{e}}

\numberwithin{equation}{section}
\begin{document}

\begin{center}
{\bf{\large Dilation of a class of quantum
dynamical
semigroups with unbounded generator on UHF algebras}}\\

{\large Debashish Goswami, Lingaraj Sahu {\footnote {The  author would like to
acknowledge the support of National Board of Higher Mathematics, DAE, India and to the
DST-DAAD programme.}}}\\

and\\
{\large  Kalyan B. Sinha {\footnote {The  author would like to
acknowledge the support from Indo-French Centre for the Promotion  of
Advanced Research as well as from DST-DAAD programme. }}}
\\
 \emph{Stat-Math Unit, Indian Statistical Institute,\\ 203,
B.T. Road, Kolkata 700 108, India.
\\ email
:  goswamid@isical.ac.in,
lingaraj\_r@isical.ac.in and kbs@isical.ac.in}
\vspace{1cm}
 \bf{ Dedicated to the memory of Professor Paul Andr{\'e} Meyer}
\end{center}
\begin{abstract}
Evans-Hudson flows are constructed for a class of quantum dynamical semigroups
with unbounded generator  on  UHF algebras, which appeared  in  \ \cite{Ma}.
It is  shown that these  flows are  unital and covariant. Ergodicity  of the  flows
for  the semigroups associated with partial states is  also discussed.
\end{abstract}

\section{Introduction}
Quantum dynamical semigroups (to be abbreviated as QDS)
 constitute a natural generalization of
classical Markov semigroups arising as expectation semigroups of  Markov
processes. A QDS $\{T_t:t\ge 0\}$ on a $C^*$-algebra $\A$ is  a $C_0$-semigroup
of completely positive   (CP) maps $T_t$ on $\A.$ Given such  a QDS, it is
interesting and  important to look for a dilation in the sense of Evan-Hudson (EH)
i.e. a family of
$*$-homomorphism  $ \eta_t:\mathcal A \rightarrow A\otimes(\Gamma
(L^2(\MBBR_+,\MBFK_0)))$  where $\MBFK_0 $ is   some separable Hilbert
space  and $\Gamma (\cdot)$ denotes the symmetric Fock
space, satisfying a suitable quantum stochastic differential equation (QSDE).
This problem has  been completely solved for  QDS with bounded generators
by Goswami, Sinha and Pal ~\cite{GS1,GPS}, where a canonical EH dilation for an
arbitrary QDS with bounded generator  has been constructed. However,
only partial success has been achieved for  QDS with unbounded generator.
It is perhaps too much to expect a complete general theory for an arbitrary QDS.
It may be wiser to look for  EH dilation for special classes of QDS. In  \cite{GS2}
for example, the author  gave a  general theory  for QDS on a $C^*$-algebra $\A,$
which is covariant with respect to an action of a  Lie group and also
symmetric with respect to  a given trace. However,  in the present article , we shall
try to construct EH dilation for another class of QDS on UHF $C^*$-algebra,
studied by T. Matsui \ \cite{Ma}. This construction has some similarity with
the earlier one but the action of the discrete group $\MBBZ^d$
instead of a Lie group  action in \  \cite{GS2}, makes the present model
somewhat  different
from that of \cite{GS2}.
We have not only proved the existence of the dilation  (in section 3),
 we  are also able to
prove that the EH flow is indeed  covariant with respect to
the $\MBBZ^d$ action (in section 4). Some ergodicity properties  of  flows are briefly
discussed too( in section 5).


\section{Notation and preliminaries}
Matsui (in  \ \cite{Ma}) constructed a class of conservative completely positive
 semigroups on
  the UHF $C^*$-algebra $\mathcal A $  generated by  infinite tensor
 product of finite
dimensional matrix algebras $ M_N(\MBBC) ,$ i.e.\ the
$C^*$-completion of $\otimes_{j \in \MBBZ^d} ~ M_N(\MBBC) ,$\
where $ N $ and $ d $ be two fixed positive integers
(inductive limit of full matrix algebras $\{ M_{N^n}(\MBBC),n\ge 1
\} $ with respect to the imbedding of $M_{N^n} $ in $ M_{N^{n+1}} $  by sending
 $ a$ to $ a \otimes 1
).$  The unique normalized trace $tr $ on $ \mathcal A $
is given  by $ tr(x)=\frac{1}{N^n} ~{Tr(x)},$ for $ x \in
M_{N^n}(\MBBC),$ where $ Tr $ denote the ordinary trace on $
M_{N^n}(\MBBC).$  For $x\in M_N(\MBBC) $ and  $ j\in \MBBZ^d, $
define an element $ x^{(j)} \in \mathcal A $ whose $ j$-th
component is $ x $ and rest are identity of $  M_N(\MBBC).$
For a simple tensor element $ a \in \mathcal A,$ let $ a_{(j)} $ be the
$j$-th component of $a$, the support of $a,$ denoted by $supp(a)$ be the set
$\{j \in
\MBBZ^d :a_{(j)}\ne 1 \}$  and for a general  element $a\in
\mathcal A ,a=\sum_{n=1}^\infty c_n a_n $ with $a_n$'s  simple
tensor elements in $\mathcal A $ and $c_n $'s  complex coefficients,
define  $supp(a)=\bigcup_{n\in \mathbb{N}}supp(a_n) $  and set $
|a|=$ cardinality of $ supp(a).$
For any $\Lambda \subseteq \MBBZ^d,\ \mathcal
A_\Lambda$ denote  the  $*$-subalgebra generated by elements of
$\mathcal A$ with support $\Lambda. $ When $\Lambda=\{k\},$ we write
$\A_k$ instead of $\A_{\{k\}}.$ Let $ {\mathcal A}_{loc} $ be
the $ *$-subalgebra of $ \mathcal A $ generated by elements $a\in
\mathcal A $ of finite support  or equivalently by $
\{x^{(j)}:x\in M_N(\MBBC),j\in \MBBZ^d\}. $ Clearly $ {\mathcal
A}_{loc} $ is dense in  $ \mathcal A. $
 \ For $k\in \MBBZ^d $ the
translation $ \tau_k $ on $\mathcal A $ is an automorphism
determined by $\tau_k (x^{(j)})=x^{(j+k)},\forall x\in M_N (\MBBC)$
and $ j\in \MBBZ^d.$ \ Thus, we get an action $ \tau $ of the infinite
discrete group $ \MBBZ^d $ on $ \mathcal A.$ \ For $x\in \mathcal
A $ denote $\tau_k(x)$ by \ $ x_k.$  The algebra $ \mathcal A $
is naturally sitting inside $ \MBFH_0 ={L^2}(\mathcal A,tr),$
the GNS Hilbert space for $({\mathcal A}, {tr}).$ It is easy to
see that $ \tau_k $ extends to a unitary on $  \MBFH_0 ,$
to be denoted  by same symbol $ \tau_k,$ giving  rise to a
unitary representation $ \tau $ of the group $ \MBBZ^d $ on
$\MBFH_0,$ which implements the action $\tau.$

We also need another dense subset of $\mathcal A,$  in a sense
like the first Sobolev space  in $ \mathcal A.$  For this, we
need to note that $M_N(\MBBC) $ is spanned by   a pair of
noncommutative representatives $ \{ U, V  \}$ of $ \MBBZ_N =\{0,1\cdots N-1\}$ such
that $ U^N=V^N=1\in M_N(\MBBC) $  and $ UV= wVU,$ where $w\in
\MBBC $ is the primitive $ N$-th root of unity( these $ U,V $ are
 given by $N\times N $ circulant matrices, note that  for $N=2,$
 $ U $ and $ V $ are respectively the Pauli-spin matrices $ \sigma_x $
 and $ \sigma_z).$
 For $ j\in\MBBZ^d $ and $\alpha,\ \beta \in G
\equiv\MBBZ_N\times\MBBZ_N,$ \ set $ \sigma_{j;\alpha,\beta}(x) = \left[
U^{(j)}V^{(j)},x\right],\ \forall x\in \mathcal A$  and $ \| x
\|_1=\sum_{j;\alpha,\beta}\|\sigma_{j;\alpha,\beta}(x) \|. $ \ Set
${\mathcal C}^1(\mathcal A)= \{x\in
\mathcal A $:\ $ {\| x \|}_1  < \infty \}.$ It is easy to see that
$\| x^* \|_1=\|\tau_j (x) \|_1=\| x \|_1$ and since
${\mathcal C}^1(\mathcal A)  $ contains the dense  $*$-subalgebra $
{\mathcal A}_{loc},\  {\mathcal C}^1(\mathcal A) $ is a  dense $\tau $
invariant
$*$-subalgebra of $
\mathcal A $. Let $ \mathcal G =\prod_{j \in \MBBZ^d} G $ be the
infinite direct product of the finite group $ G $ at each lattice
site. Thus each $ g\in \mathcal G $ has $j$-th component $g_{(j)}
= (\alpha_j,\beta_j)$ with $ \alpha_j,\beta_j \in G $  and  for $g
\in \mathcal G $ define its support by $ supp(g) = \{j\in
\MBBZ^d:g_{(j)}\neq(0,0)\} $ and $ \left|g\right|=$  cardinality
of $ supp(g).$   Consider the projective unitary representation of
$ {\mathcal G},$ given by $ {\G}\ni g {\mapsto} {U_g } = \prod_{j \in
\MBBZ^d}{U^{(j)}}^{\alpha_j}{V^{(j)}} ^{\beta_j} \in \mathcal A .$

For a given CP map $T$ on $\A,$ formally we  define the Linbladian \\
\[{\mathcal
L}=\sum_{k\in \MBBZ^d }{\mathcal L}_k,\]
\[ \mbox{ where } \ {\mathcal
L}_k x=\tau_k \mathcal L_0(\tau_{-k}x ),\ \forall x\in \A_{loc}\]

\be \label{linblad} \mbox{with}\ \mathcal L_0(x)=-\frac{1}{2}\{T(1),x\}+T(x),\ee
where
$\{A,B\}:=AB+BA.$\\
In particular consider the Linbladian  $\mathcal L$ for the CP map
\[Tx=\sum_{l=0}^\infty a_l^*xa_l, \ \forall x\in \A,\]associated with
  a sequence
of elements $\{a_l\} $ in $\A,\ a_l=\sum_{g\in\G} c_{l,g}U_g $ such
that \\
$ \sum_{l=0}^\infty \sum_{g\in\G}\left| c_{l,g}\right|{\left| g \right|}^2 <
\infty.  $  Matsui has proven the following in the
paper referred  earlier \ \cite{Ma} .
\bthm\label{semigroup}
(i)The $ \mathcal L $ formally defined above
 is well defined on $
\C^1(\A)$ and the closure $\hat{\mathcal L} $ of $ {\mathcal
L}/_{\C^1(\A)} $ is a generator of conservative CP semigroup $
\{P_t:t\ge 0\}$ on $ \A,$\\
(ii)\ The semigroup $ \{P_t\} $ leaves  $ \C^1(\A) $ invariant.
\ethm
\noindent The  semigroup $ P_t $ satisfies
 \[P_t(x)=x+\int_0^tP_s  ( \hat{\mathcal L} (x))ds,\
 \forall x\in Dom(\hat{\mathcal L}).\]
Since  $1 \in \C^1(\A) $ ( in fact $\| 1 \|_1=0 $ ) and
$\hat{\mathcal L}(1)=\mathcal L(1)=0,$ it follows that $ P_t(1)=1,\forall t\ge0.$

Following  \cite{Ma}, we say that $P_t$ is ergodic if
there exist an invariant
state $\psi$ satisfying
\be \label{ergodic}\| P_t(x)-\psi(x)1\|\rightarrow 0 \
\mbox{as}\ t\rightarrow \infty,\ \forall x\in {\mathcal A}.\ee
In  \ \cite{Ma},\ the  author has discussed some criteria for ergodicity of
 CP semigroup  $P_t.$
Some examples of such semigroups, associated with partial
states  on the UHF algebra and their perturbation are given.

Let $\phi $ be a state on $ M_N(\MBBC)$ and for  $k\in \MBBZ^d, $ the partial
state $\phi_k$ on $\mathcal A $ determined by
$\phi_k(x)=\phi(x_{(k)})x_{\{k\}^c},$  for $x=x_{(k)} x_{\{k\}^c},$ where
$x_{(k)} \in \mathcal A_k $ and $ \ x_{\{k\}^c}\in \mathcal
A_{\{k\}^c}.$ We can find  elements $\{L^{(m)}:m=1,2\cdots N^\prime\}$
 in $M_N(\MBBC),$ for some finite natural number $N^\prime$ such that
\[\phi(x)=\sum_{m=1}^{N^\prime} { L^{(m)}}^*xL^{(m)}\ \forall x\in M_N(\MBBC) \ \mbox{and}
 \sum_{m=1}^{N^\prime} { L^{(m)}}^*L^{(m)}=1\]
 For  $ m=1,\cdots N^\prime, $
 consider the element $ L_0^{(m)}\in \A_0 $
 with zeroth component  is $  L^{(m)}$ respectively. Now for $k\in \MBBZ^d$
 and $ m=1,\cdots N^\prime,$
 writing $ L_k^{(m)}= \tau_k( L_0^{(m)}),$ the partial state  $\phi_k$ is given by,
\[\phi_k(x)=\sum_{m=1}^{N^\prime} { L_k^{(m)}}^*xL_k^{(m)}\ \forall x\in \A\]
By (\ref{linblad}), formally the Linbladian $\mathcal L^\phi, $
corresponding to the partial state $\phi_0 $ is given  by
\[\mathcal L^\phi(x)= \sum_{k\in
 \MBBZ^d}\mathcal L_k^\phi(x)\]
 where \[\mathcal L_k^\phi(x)=\phi_k(x)-x=
 \frac{1}{2}\sum_{m=1}^{N^\prime} [{ L_k^{(m)}}^*,x]L_k^{(m)}+
 { L_k^{(m)}}^*[x,L_k^{(m)}]\]
 It  follows from theorem (\ref{semigroup})  that
$\mathcal L^\phi$  defined on  $\C^1(\A).$
 Moreover the closure $\hat{\mathcal
L}^\phi$ of $ {\mathcal
L^\phi}/_{\C^1(\A)} $  generates a conservative
 CP semigroup $ P_t^\phi$ on $\mathcal A$
given by
\[P_t^\phi(\prod_{k\in \Lambda}x_{(k)})=\prod_{k\in \Lambda}\
\{\phi(x_{(k)})+e^{-t}(x_{(k)}-\phi(x_{(k)}))\}.\]
Note that the map $ \Phi $ define by,
\[ \Phi(\prod_{k\in \Lambda}x_{(k)})=\lim_{t\rightarrow \infty}
 P_t^\phi(\prod_{k\in \Lambda}x_{(k)})=\prod_{k\in \Lambda}\phi(x_{(k)})\]
 extends as a state on $\mathcal A$  which is  the unique ergodic
  state  for $ P_t^\phi.$
For any real $c,$ consider the perturbation

\[\mathcal
L^{(c)}(x)=\mathcal
L^\phi(x)+c\mathcal L(x), \forall x\in \C^1(\A).\]
It clear that $L^{(c)}$ is the Linbladian associated with the CP map
\[T(x)=\sum_{m=1}^{N^\prime} { L_k^{(m)}}^*xL_k^{(m)}+c \sum_{l=0}^\infty a_l^*xa_l,
 \forall
x\in \A\]
and by theorem (\ref{semigroup})
it follows that the closure
$\hat{\mathcal L^{(c)}}$ of $\mathcal
L^{(c)}/\C^1(\A)$ generate a QDS $P_t^{(c)}.$ Moreover,\\
\bthm \cite{Ma}
There exist a constant $c_0$ such that for $0\le c\le c_0$ the above semigroup
$P_t^{(c)}$ have the  unique ergodic state $\Phi^{(c)}$ and
\be \label{perturb}\|P_t^{(c)}(x)\|_1\le 2e^{-(1-\frac{c}{c_0})t}\|x\|_1\ \mbox{and}\ee
\[\|P_t^{(c)}(x)-\Phi^{(c)}(x)1\|\le \frac{4}{N^2}e^{-(1-\frac{c}{c_0})t}\|x\|_1,\
 \forall x\in \C^1(\A).\]
\ethm
\brmrk
The ergodic state $\Phi^{(c)}$ corresponding  QDS $P_t^{(c)}$ is given by
\[\Phi^{(c)}(x)=\Phi(x)+ c \int_0^\infty\Phi({\mathcal L}( P_t^{(c)}(x))) dt,
\forall x\in \C^1(\mathcal A)\]

\ermrk

Let us conclude the present section with a brief discussion on  the
 fundamental integrator process of
quantum stochastic calculus, introduced by Hudson and Parthasarathy \cite{HP}).
Let $ \MBFK=L^2(\MBBR_+,\MBFK_0)$ where   $ \MBFK_0 =l^2(\MBBZ^d)$ with
 the canonical orthonormal  basis $\{e_j :j \in \MBBZ^d \}$
and $\Gamma =\Gamma_{sym}(\MBFK),$
 the symmetric Fock space over $\MBFK.$ For $f\in \MBFK,$ we
 denote by $ \EXP(f) $ the exponential vector in $ \Gamma $ associated
 with  $ u \ :$\\
 \[\EXP(f)=\bigoplus_{n\ge 0}\frac{1}{\sqrt{n!}}u^{(n)}\]
 where $ u^{(n)}=\underbrace{u \otimes u  \otimes
 \cdots \otimes  u}_{n-copies} $ for $ n>0 $
and by convention  $u^{(0)}=1 .$
 For $ f=0,\  \EXP(f) $ is called the vacuum vector in $\Gamma.$
 Let $\mathcal C $ be the space of all bounded
continuous functions on $L^2(\MBBR_+,\MBFK_0),$ so that  $\mathcal
E(\mathcal C)\equiv\{ \EXP(f):f\in \mathcal C \}$ is  total in
$\Gamma(\MBFK).$  Any  $ f\in L^2(\MBBR_+,\MBFK_0)$ decomposes as
 $ f=\sum_{k\in \MBBZ^d}f_k e_k $  with $ f_k \in L^2(\MBBR_+).$
We take the freedom to use the same symbol
$ f_k $ to denote  function in $ L^2(\MBBR_+,\MBFK_0)$ as well, whenever
  it is clear from the context. The  family
  of fundamental processes,
 $ \{\Lambda_i^j:i,j\in
 \MBBZ^d\},$ associated with the orthonomal basis $ \{e_j
:j \in \MBBZ^d \},$ given by
 \bean \lefteqn{ \Lambda_j^i(t)}
&&~~=a_{\chi_{[0,t]}\otimes e_i} \ \mbox{for} \  i\ne 0,j=0 \
\mbox{(annihilation)}\\
&&~~=a^\dag_{\chi_{[0,t]}\otimes e_j} \ \mbox{for} \ i=o, j\ne 0\
\mbox{(creation)}\\
&&~~=\Lambda_{M_{\chi_{[0,t]}}\otimes |e_j><e_i|} \ \mbox{for} \
i,j\ne 0 \ \mbox{(conservation )}\\
&&~~= t1 \ \mbox{for} \  i=j=0,\ \mbox{(time )}\eean where $
M_{\chi_{[0,t]}}$ is the multiplication operator   on $
L^2(\MBBR_+)$ by characteristic function of the interval $ [0,t].$
For detail see \cite{KRP,Mey}).

\section{Evans-Hudson (EH)type dilation}
Formally, we would like to solve the following quantum
stochastic differential equation (QSDE) in \[
\mathcal B(L^2(\mathcal A, tr))\otimes
\mathcal B(\Gamma(L^2(\MBBR_+,\MBFK_0))),\]
 \be \label{EH}
dj_t(x)=\sum_{j\in \MBBZ^d}j_t(\delta _j^\dag (x))
d{a_j}(t)+\sum_{j\in \MBBZ^d}j_t(\delta_j (x))
da_j^\dag(t)+j_t({\hat{\mathcal L}} (x))dt, \ee
\[{j_0}(x)=x \otimes 1_{\Gamma}\ .\]
Now if  we look at the corresponding  (see \ \cite{KRP,Mey})\
 Hudson-Parthasarathy (HP) equation
in $ L^2(\mathcal A, tr)\otimes
L^2(\Gamma(L^2(\MBBR_+,\MBFK_0))),$
 \be \label{HPeq}
dU_t=\sum_{j\in \MBBZ^d}\{r_j^*
d{a_j}(t)-r_j
da_j^\dag(t)-\frac{1}{2}r_j^*r_j dt\}U_t,\ee
\[{U_0}(x)=1_{L^2 \otimes \Gamma}.\]
However, though each $r_j\in \mathcal A $ and hence is in
$ \mathcal B(L^2(\mathcal A, tr)), $ the equation \ (\ref{HPeq})
does not admit a solution since
\[\langle u,\sum_{j\in \MBBZ^d}r_j^*r_j u\rangle=\sum_{j\in \MBBZ^d}\|r_j u\|^2\]
is not convergent in general and hence $ \sum_{j\in \MBBZ^d}r_j\otimes e_j $ does not
define an element in $ \mathcal A \otimes \MBFK_0.$ For example, let $r $
be the singleton $ U^{(k)}\in \mathcal A $ so that $ r_j=U^{(k+j)}$ is a unitary
for all $j\in \MBBZ^d$
and hence
\[\sum_{j\in \MBBZ^d}\|r_j u\|^2=\sum_{j\in \MBBZ^d}\|u\|^2 =\infty. \]
However, as we shall see, in many situation
there exist EH flows, even though the corresponding HP equation \ (\ref{HPeq})\
does not admit  a solution.

\brmrk
 There are some cases when an EH dilation can be seen to be implemented
by  solution of HP equation, for example, given  a self adjoint $r\in \A,$
\[dV_t=\sum_{k\in
\MBBZ^d}V_t ( S_k^*
d{a_k}(t)-S_k
da_k^\dag(t)-\frac{1}{2}S_k^*S_k dt),\  V_0=1,\]
where $S_k$ is  defined by $S_k(x)=[r_k,x]$ for $x\in \A
\subseteq L^2(\mathcal A, tr),$ admits a unique unitary solution  and
\[x\mapsto V_t^*(x\otimes 1)V_t\] gives an EH dilation for $P_t$ (ref  \cite{Mo,MS}).
\ermrk
Let $a,b\in\MBBZ_N$ be fixed, $ W={U^a}{V^b}\in \M_N(\MBBC)$ a
fixed element. Consider the following
representation of the infinite product  group
$ {\mathcal G}^\prime=\prod_{j \in \MBBZ^d}\MBBZ_N,$  given by
$${\mathcal G}^\prime\ni g {\mapsto}{W_g }= \prod_{j \in
\MBBZ^d}{W^{(j)}}^{\alpha_j} ,\ \mbox{where}\  g=({\alpha_j}).$$\\
For any $y \in \A,\ y= \sum_{g\in \G}c_gU_g $ and for
$n\ge 1$ we define
\[\vartheta_n(y)=\sum_{g \in \G} \left|c_g\right|{\left|g\right|}^n. \]
Now consider $ r \in \A,\ r= \sum_{g\in {\G}^\prime}c_gW_g\ $ such that
$ \sum_{g\in\G^\prime}\left| c_{g}\right|{\left| g \right|}^2 <
\infty.$ It is clear that $ \vartheta_1(y)=\sum_{g \in \G^\prime} \left|c_g\right|
{\left|g\right|}< \infty. $  Note that for  $ x\in {\mathcal A}_{loc} $ we can always
write  $ x=\sum_{h\in \G }c_hU_h,$  with complex coefficients $c_h$ satisfying
$c_h=0 $ for  $ \forall h $ such  that $ supp(h)\bigcap supp(x)$ is  empty. So
\[\vartheta_n(x)=\sum_{h \in \G} \left|c_h\right|{\left|h\right|}^n < \infty
\mbox{ for
$n\ge 1$}\]
 and it is clear that
\[\vartheta_n(x)\le |x|^n\sum_{h \in \G} |c_h|\le c_x^n \]
 for constant $ c_x =|x|(1+\sum_{h \in \G} |c_h|).$
Now consider the formal Linbladian define by (\ref{linblad})
$\mathcal L$ associated with
CP map $T(x)=r^*xr,$
 ${\mathcal
L}=\sum_{k\in \MBBZ^d }{\mathcal L}_k,$ where $ {\mathcal
L}_k x=\tau_k \mathcal L_0(\tau_{-k}x )$ with
 ${\mathcal
L}_{0}x=\frac{1}{2} \{
\left[r^*,x\right]r+r^*\left[x,r\right]\}$ so that,
 \[ \mathcal L_k(x)=\frac{1}{2}\{\left[r_k^*,x\right]r_k+r_k^*\left[x,r_k\right]\}.\]
  Let denote  these two bounded  derivations $ \left[r_k^*,
.\right] $ and $ \left[ .,r_k\right] $  in $ \A,$  \    by $
\delta_k^\dag $ and $ \delta_k $ respectively so that ${ \mathcal
L}(x)= \frac{1}{2}\sum_{k \in \MBBZ^d}\delta_k^\dag(x)r_k+
r_k^*\delta_k(x).$

 For $ n\ge
1,$  denote the set of integers $\{1,2,\cdots n\} $ by $
I_n$ and for  $ 1\leq p\leq n,\ P=\{l_1,l_2 \cdots l_p\}\subseteq
I_n $ (where $ l_i $'s are in  increasing order), define a map from
$n$-fold Cartesian product of  $ \MBBZ^d$ to that  $ p $ copies of $ \MBBZ^d$
given by
\[ \bar{k}(I_n)=(k_1,k_2 \cdots k_n)\mapsto
\bar{k}(P)=(k_{l_1},k_{l_2} \cdots k_{l_p})\] and similarly,
$\bar{\varepsilon}(P)=(\varepsilon_{l_1}\varepsilon_{l_2}
\cdots\varepsilon_{l_p})$ for  a vector
$\bar{\varepsilon}(I_n)=(\varepsilon_1\varepsilon_2
\cdots\varepsilon_n)$ in  $n$-fold Cartesian product of $\{-1,0,1\}.$\\
For brevity of notations, we write
 $ \bar{\varepsilon}(P)\equiv c\  (c\in \{-1,0,1\})
 $ to  mean that
all  $ \varepsilon_{l_i}=c$ and $ \bar{k}(n),\bar{\varepsilon}(n)$ will stand for
 $\bar{k}(I_n),\bar{\varepsilon}(I_n)$
respectively. Setting  $
\delta_k^\varepsilon =\delta_k^\dag,{\mathcal L}_k
 $ and $\delta_k $ depending upon $ \varepsilon =-1,0 $\  and   $1
$ respectively, we write $R(\bar
{k})=r_{k_1}r_{k_2}\cdots r_{k_p}$  and
$\delta(\bar{k},\bar{\varepsilon})= \delta^{\varepsilon_p}_{k_p} \cdots
{\delta^{\varepsilon_1}_{k_1}}$  for any $\bar
{k}=(k_1,k_2\cdots k_p)$ and  $\bar
{\varepsilon}=(\varepsilon_1,\varepsilon_2\cdots \varepsilon_p).$
Now we have the following useful lemma,
\blema\label{lema1}Let $r, x $ and constant  $c_x$
be as above.\ Then
\\(i)For any $n\ge1,$ \[\sum_{\bar{k}(n)}\|
\delta(\bar{k}(n),\bar{\varepsilon}(n))(x)\| \leq ({2\vartheta_1(r)
c_{x}})^n, \forall x\in {\A}_{loc}, \]
where  $ \bar{\varepsilon}(n) $ is  such that
$\varepsilon_l\neq 0,\ \forall \ l \in I_n.$
\\(ii) For any $n\ge 1 $ and $  \bar{k}(n), $
 \bean \lefteqn
{ {\mathcal L}_{k_n}\cdots{\mathcal L}_{k_1}(x)}\\
&=&\frac{1}{2^n} \sum_{p=0,1\cdots n}\sum_{P\subseteq
I_n:|P|=p}{{R(\bar
{k}({P^c}))}^*}\delta(\bar{k}(n),{\bar{\varepsilon}}_{(P)}(n))(x)R(\bar
{k}(P)),\\
\eean
 where ${\bar{\varepsilon}}_{(P)}(n)$ is such that $
{\bar{\varepsilon}}_{(P)}(P)\equiv -1 $ and $
{\bar{\varepsilon}}_{(P)}({P^c})\equiv 1.  $
\\(iii)For any $ n \geq 1,p \leq n, P \subseteq
I_n $ and $ \bar{\varepsilon}(n) $ is such that $ \bar{\varepsilon}(P) $
contains  all  those components having value $ 0,$ we have,
 \[\sum_{\bar{k}(n)}\|\delta(\bar{k}(n),\bar{\varepsilon}(n))(x)\| \leq {\|r\|}^{p}
({2\vartheta_1(r) c_{x}})^n \]
\[ \leq {(1+\|r\|)}^n ({2\vartheta_1(r) c_{x}})^n. \]
\\(iv) Let $ m_1,m_2 \geq 1;\ x,y\in \mathcal A_{loc} $  and
$ {\bar{\varepsilon}^\prime}(m_1),
{\bar{\varepsilon}}^{\prime\prime}(m_2)  $ be two fixed tuples,
then for $n\ge 1$  and $ \bar{\varepsilon}(n) $ as in (iii),\ we have,
\[\sum_{\bar{k}(n),{\bar{k}}^\prime
(m_1),{\bar{k}}^{\prime\prime}(m_2)} \|
\delta(\bar{k}(n),\bar{\varepsilon}(n))\{
\delta({\bar{k}}^\prime(m_1),
{\bar{\varepsilon}}^\prime(m_1))(x) \cdot
\delta({\bar{k}}^{\prime\prime}(m_2),{\bar{\varepsilon}}^{\prime\prime}(m_2))(y)
\} \|\]
\[ \leq {2^n}{(1+\|r\|)}^{2n+m_1+m_{2}} ({2\vartheta_1(r)
c_{x,y}})^{n+m_1+m_{2}},\] \ where $ c_{x,y}= max
\{c_{x},c_{y}\}.$
\elema
\begin{proof}
(i) As  $r^*$ is again of the same  form as $r,$ it is enough to
observe the following
\[\sum_{k_n,\cdots k_1}\| \left[r_{k_n},\cdots
\left[r_{k_1},x\right]\cdots \right]\|  \le ({2\vartheta_1(r)
c_{x}})^n , \forall x\in {\A}_{loc} \ . \]
In order to prove this  consider
\[ LHS=\sum_{k_n,\cdots k_1}\ \sum_{{g}_n,\cdots{g}_1 \in {\G}^\prime;
h\in{\G}} |c_{g_n}| \cdots |c_{g_1}|\  |c_h|\
\|\left[{\tau_{k_n}}{W_{g_n}},\cdots
\left[{\tau_{k_1}}{W_{g_1}},U_h\right]\cdots\right]\|.\]

\noindent  Since for any two commuting elements $ A, B $ in $\A,\
\left[A\left[B,x\right] \right]=\left[B\left[A,x\right] \right] ,$
for the commutator
$\left[{\tau_{k_n}}{W_{g_n}},\cdots
\left[{\tau_{k_1}}{W_{g_1}},U_h\right]\cdots\right]  $ to be nonzero, it is necessary
to have,\ $ (supp(g_i)+k_i)\bigcap
supp(h) \neq \phi  $  for each  $i=1,2,\cdots n $
 and number  of choices of
such $k_i\in \MBBZ^d $ is at most $|g_i|\cdot |h|. $ Thus we get,
\bean \lefteqn{ \sum_{k_n, \cdots k_1}\|\left[r_{k_n}, \cdots
\left[r_{k_1},x\right] \cdots \right]\|}\\
&&\leq \sum_{{g}_n,\cdots{g}_1 \in {\G}^\prime;
h\in \G} |c_{g_n}|\cdots|c_{g_1}| |c_h||{g_n}|\cdots|{g_1}|
{|h|}^n2^n\\
&&\leq ({2\vartheta_1(r)
c_{x}})^n.  \eean
(ii)The proof is by induction. For any\ $ k\in \MBBZ^d $\ we have,
\[{{\mathcal L}_k}(x)= \frac{1}{2}\sum_{k \in\MBBZ^{d}}\delta_k^\dag(x)r_k+
r_k^*\delta_k(x),\] so it is trivially  true for $ n=1 .$\
 \ Assume it to be true for some  $ m > 1 $  and  for  any $
k_{{m}+1}\in \MBBZ^d $ consider $ {\mathcal
L}_{k_{{m}+1}}{\mathcal L}_{k_m}\cdots{\mathcal L}_{k_1}(x),$  \
now by applying the statement for $n= m $  we get,
 \bean\lefteqn { {\mathcal L}_{k_{m+1}}{\mathcal L}_{k_m}\cdots{\mathcal
L}_{k_1}(x)}\\
&=& \frac{1}{2^{m+1}} \sum_{p=0,1\cdots m}\ \sum_{P\subseteq
I_m:|P|=p}[\delta_{k_{m+1}}^*\{{{R(\bar
{k}({P^c}))}^*}\delta(\bar{k}(m),{\bar{\varepsilon}}_{(P)}(m))(x)R(\bar
{k}(P))\}r_{k_{m+1}} \\
&&~~~~~~~~~~~~~~~~~~~~~~~~~~~~~+  r_{k_{m+1}}^*\delta_{k_{m+1}}\{{{R(\bar
{k}({P^c}))}^*}\delta(\bar{k}(m),{\bar{\varepsilon}}_{(P)}(m))(x)R(\bar
{k}(P)) \}],\\
\eean
 Since  $r_k$'s are commuting  with each other,
\bean \lefteqn{LHS}\\&&=\frac{1}{2^{m+1}} \sum_{p=0,1\cdots m}\ \sum_{P\subseteq
I_m:|P|=p}[{{R(\bar
{k}({P^c}))}^*}\delta_{k_{m+1}}^*
\delta(\bar{k}(m),{\bar{\varepsilon}}_{(P)}(m))(x)R(\bar
{k}(P))r_{k_{m+1}}\\
&&~~~~~~~~~~~~~~~~~~~~~~~~~~~~~~~~+ r_{k_{m+1}}^* {{R(\bar
{k}({P^c}))}^*}\delta_{k_{m+1}}\delta(\bar{k}(m),{\bar{\varepsilon}}_{(P)}(m))(x)R(\bar
{k}(P))] \\
&&=\frac{1}{2^{m+1}} \sum_{p=0,1\cdots m+1}\ \sum_{P\subseteq
I_{m+1}:|P|=p}{{R(\bar
{k}({P^c}))}^*}\delta(\bar{k}(m+1),{\bar{\varepsilon}}_{(P)}(m+1))(x)R(\bar
{k}(P)). \eean
(iii) By simple application of (ii),
\[\delta(\bar{k}(n),\bar{\varepsilon}(n))(x)\]
\be \label{LL} =\frac{1}{2^p} \sum_{q=0,1\cdots p}\ \sum_{Q \subseteq
P:|Q|=q}{{R(\bar {k}({P\setminus
Q}))}^*}\delta(\bar{k}(n),{\bar{\varepsilon}}_{(Q,P)}(n))(x)R(\bar
{k}(Q)), \ee
 where ${{\bar{\varepsilon}}}_{(Q,P)}(n) $ define by the map
 from the  $n$-fold Cartesian product of $\{-1,0,1\}$ to
 itself,
 $ \bar{\varepsilon}(n)\mapsto {{\bar{\varepsilon}}_{(Q,P)}}(n)$ such that
 ${{\bar{\varepsilon}}_{(Q,P)}}(Q)\equiv  -1,
{{\bar{\varepsilon}}_{(Q,P)}}({P\setminus Q})\equiv 1$ and
\\ ${{\bar{\varepsilon}}_{(Q,P)}}(I_n\setminus P)=
\bar{\varepsilon}(I_n\setminus P),$
now (i) gives  that we  require .
\\(iv) By \ (\ref{LL}) we have,
\bean \lefteqn{ LHS}\\ &&=\frac{1}{2^p}\sum_{\bar{k}(n),{\bar{k}}^\prime
(m_1),{\bar{k}}^{\prime\prime}(m_2)}\ \sum_{q=0,1\cdots
p}\ \sum_{Q \subseteq P:|Q|=q} \| R(\bar {k}({P\setminus Q}))^*\\
&&~~~~~\delta(\bar{k}(n),{\bar{\varepsilon}}_{(Q,P)}(n))\
[ \delta({\bar{k}}^\prime(m_1),
{\bar{\varepsilon}}^\prime(m_1))(x)\cdot
\delta({\bar{k}}^{\prime\prime}(m_2),{\bar{\varepsilon}}^{\prime\prime}(m_2))(y)
]\  R(\bar {k}(Q))\|,\eean Now applying  Leibnitz rule,  it become
\bean \lefteqn{\leq \frac{{\|r\|}^{p}}{2^p}\sum_{\bar{k}(n),{\bar{k}}^\prime
(m_1),{\bar{k}}^{\prime\prime}(m_2)}\ \sum_{q=0,1\cdots
p}\ \sum_{Q \subseteq P:|Q|=q}\ \sum_{l=0,1\cdots n}\ \sum_{L \subseteq
I_n:|L|=l}}\\
&&\|\delta(\bar{k}(L),{\bar{\varepsilon}}_{(Q,P)}(L))
\delta({\bar{k}}^\prime(m_1),
{\bar{\varepsilon}}^\prime(m_1))(x)\|\\
&&\|\delta(\bar{k}(L^c),
{\bar{\varepsilon}}_{(Q,P)}(L^{c}))
[\delta({\bar{k}}^{\prime\prime}(m_2),{\bar{\varepsilon}}^{\prime\prime}(m_2))(y)]\|.
\eean

Now by  using (iii), we  obtain,
\bean \lefteqn  { LHS}\\
&& \leq \frac{{(1+\|r\|)}^n}{2^p} \sum_{q=0,1\cdots
p} c(p,q) \sum_{l=0,1\cdots n} c(n,l) {(1+\|r\|)}^{l+m_1}
({2\vartheta_1(r) c_{x}})^{l+m_1}\\
&&~~~~~~~~~~~~~~~~~~ \cdot {(1+\|r\|)}^{n-l+m_2}
({2\vartheta_1(r) c_{y}})^{n-l+m_2}\ \mbox{(where $c(p,q)=\frac{p!}{(p-q)!q!}$)}\\
&& ~~~~~~~~~~~~~~~\leq {2^n}{(1+\|r\|)}^{2n+m_1+m_{2}} ({2\vartheta_1(r)
c_{x,y}})^{n+m_1+m_{2}}. \eean
\end{proof}
\noindent Now we are in position to prove the following result about existence
 of EH dilation of CP semigroup $P_t$ associated with element $r\in \A$ discussed above.
\bthm \label{theorem1} For $t\geq 0 $ and $x,y \in {\mathcal A}_{loc},$ \\
(a)\ There exist unique solution $j_t(x) $ of the QSDE,  \be \label{qsde}
dj_t(x)=\sum_{j\in \MBBZ^d}j_t(\delta_j^\dag x)
d{a_j}(t)+\sum_{j\in \MBBZ^d}j_t(\delta_j x)
d{{a^\dag}_j}(t)+j_t({\hat{\mathcal L}} x)dt, \ee
\[{j_0}(x)=x \otimes 1_{\Gamma}\ .\]
such that $ j_t(1)=1,\ \forall t\ge 0. $
\\(b)\ For $u,v \in\MBFH_0
,f,g\in \mathcal C ,$ \ \be \label{homo1} \langle
u\EXP(f),j_t(xy)v\EXP(g) \rangle=\langle j_t(x^{*}) u\EXP(f),j_t(y)v\EXP(g)
\rangle.  \ee
\\(c)\ $j_t $ is a contraction on $ {\mathcal A}_{loc}
$ and  extends uniquely to a unital $C^*$-homomorphism from
$\mathcal A $ in to $\mathcal A \otimes \mathcal B (\Gamma).$
\ethm
\begin {proof} Note first that $ {\mathcal A}_{loc}  $ is  a dense $*$-sub algebra
of $ {\mathcal A}. $
\\(a) As usual, we solve by iteration. For $ {t_0}\geq 0, t\leq t_0 ,
 x\in {\mathcal A}_{loc},u\in\MBFH_0
,f\in \mathcal C,$ we set
\bean \lefteqn{{j_t^{(0)}}(x)=x \otimes 1_{\Gamma}\ \mbox{ and }}\eean\\
$j_t^{(n)}(x)=x \otimes 1_\Gamma$
\be  \label{exist1}+\int_0^t \sum_{j\in \MBBZ^d}
j_s^{(n-1)}(\delta _j^\dag (x)) da_j(s)+\sum_{j\in
\MBBZ^d}j_s^{(n-1)}(\delta _j (x))
da^\dag_j(s)+j_s^{(n-1)}({\hat{\mathcal L}} (x))ds. \ee
 Then we shall show,
\[\|\{j_t^{(n)}(x)-j_t^{(n-1)}(x)\}u\EXP(f)\| \]
\be \label{exist2}\leq \frac{(t_0c_f)^{n/2}}{\sqrt{n!}}\|u\EXP(f)\|
\sum_{\bar{k}(n)}
\sum_{\bar{\varepsilon}(n)}\|\delta(\bar{k}(n),\bar{\varepsilon}(n))
(x)\| ,\ee
where $c_f=2e^{\gamma_{f}({t_0})}(1+{\|f\|}_\infty^2),$
\ with ${\gamma_{f}({t_0})}=\int_0^{t_0}(1+\|f(s)\|^2)ds. $
For $ n=1, $  by the basic estimate of quantum stochastic integral \
(\ see \cite{KRP,Mey}),
 \bean \lefteqn{\|\{j_t^{(1)}(x)-j_t^{(0)}(x)\}u\EXP(f)\|^2}\\
 && =\| \{\int_0^t\sum_{j\in \MBBZ^d}\delta _j^\dag (x)
d{a_j}(s)+\sum_{j\in \MBBZ^d}\delta _j (x) da^\dag_j(s)+{\hat{\mathcal L}}
(x) ds\}u\EXP(f)\|^2\\
&& \leq  2e^{\gamma_{f}({t_0})}{\|\EXP(f)\|}^2\int_0^t\{ \sum_{j\in \MBBZ^d}
\|\delta _j^\dag (x) u \|^2 +\sum_{j\in \MBBZ^d}{\|\delta _j
(x)u\|}^2 +{\|{\hat{\mathcal L}} (x) u\|}^2\}{(1+\|f(s)\|)}^2 ds\\
&&\leq c_f t_0 {\|\EXP(f)\|}^2\{ \sum_{j\in \MBBZ^d} \|\delta _j^\dag
(x) u\|+\|\delta _j (x)u\|+\|{{\mathcal L}_{j}} (x) u\|\}^2. \eean
Thus\ (\ref{exist2}) is true for $ n=1. $ Inductively  assume the estimate
for some $ m >1,$  again by same argument as above,
\bean \lefteqn{
\|\{j_t^{(m+1)}(x)-j_t^{(m)}(x)\}u\EXP(f)\|^2}\\
&= & \|\{ \int_0^t \sum_{j\in \MBBZ^d} [j_{s_m}^{(m)} (\delta
_j^\dag (x))-j_{s_m}^{(m-1)}(\delta _j^\dag (x))]
d{a_j}(s_{m})\\
&&\ + \sum_{j\in \MBBZ^d}[j_{s_m}^{(m)}(\delta _j
(x))-j_{s_m}^{(m-1)}(\delta _j (x))] d{{a^\dag}_j}(s_{m})\\
&&\ +[j_{s_m}^{(m)}({\hat{\mathcal L}} (x))-j_{s_m}^{(m-1)}
({\hat{\mathcal L}} (x))] ds_{m}\} u\EXP(f)\|^2\\
& \leq & 2e^{\gamma_{f}({t_0})}\int_0^t \{ \sum_{j\in \MBBZ^d}{\|
[j_{s_m}^{(m)}(\delta _j^\dag (x))-j_{s_m}^{(m-1)}
(\delta _j^\dag (x))]u\EXP(f)\|}^2 \\
&&\ + \sum_{j\in \MBBZ^d} \| [j_{s_m}^{(m)}(\delta _j
(x))-j_{s_m}^{(m-1)}(\delta _j (x))]u\EXP(f)\|^2\\
&&\ +{\| [j_{s_m}^{(m)}({\hat{\mathcal L}}
(x))-j_{s_m}^{(m-1)}({\hat{\mathcal L}}
(x))]u\EXP(f)\|}^2\}(1+{\|f(s_{m})\|}^2) ds_{m}\\
& \leq  & c_f{\int_0}^t {[ \sum_{j\in \MBBZ^d} \{ {\|
[{j_{s_m}^{(m)}}({\delta_j}^\dag (x))-{j_{s_m}^{(m-1)}}
(\delta _j^\dag (x))]u\EXP(f)\|}}\\
&&\ +\sum_{j\in \MBBZ^d} {\| [{j_{s_{m}}^{(m)}}({\delta_j}
(x))-{j_{s_{m}}^{(m-1)}}({\delta_{j}}(x))]u\EXP(f)\|}\\
&&\ +{
{\| [j_{s_m}^{(m)}({\hat{\mathcal L}} (x))-j_{s_m}^{(m-1)}({\hat{\mathcal L}}
(x))]u\EXP(f)\|}\} ]}^2ds_{m}. \eean
Now applying  \ (\ref{exist2})\ for $n=m,$ \ we get the require estimate for $
n=m+1 $ and
furthermore   by the estimate of lemma \ (\ref{lema1} (iii)),\

 \[\|\{j_t^{(n)}(x)-j_t^{(n-1)}(x)\}u\EXP(f)\| \leq
 3^n{(t_0c_f)^{n/2}}\|u\EXP(f)\|
{(1+\|r\|)}^n (1+{2\vartheta_1(r) c_{x}})^n ,\] thus  it follows that
the sequence $ \{ j_t^{(n)}(x)u\EXP(f)\} $ is cauchy.
 Define $ j_t(x)u\EXP(f) $\ to be the\ $ \lim_{n\rightarrow \infty}
{j_t}^{(n)}u\EXP(f),$  \    that is \be \label{exist3} j_t(x)u\EXP(f)=xu
\otimes\EXP(f)+\sum_{n\ge 1}\{j_t^{(n)}(x)-j_t^{(n-1)}(x)\}u\EXP(f) \ee
 and one has
\be\label{exist4} \|j_t(x)u\EXP(f) \| \leq \|u\EXP(f)\|[\|x\|+\sum_{n\ge
1}3^n {(t_0c_f)^{n/2}} {(1+\|r\|)}^n (1+{2\vartheta_1(r) c_{x}})^n].
\ee
Uniqueness follows by setting,
 \[ q_t(x)=j_t(x)-j_t^\prime(x)\]and observing
\[dq_t(x)=\sum_{j\in \MBBZ^d}q_t(\delta _j^\dag (x))
d{a_j}(t)+\sum_{j\in \MBBZ^d}q_t(\delta _j (x))
d{{a^\dag}_j}(t)+q_t({\mathcal L} (x))dt,\  {q_0}(x)=0. \]
Exactly similar estimate as above will show that, for all $ n\ge 1, $

\[\|q_t(x)u\EXP(f)\|
\leq \frac{(t_0c_f)^{n/2}}{\sqrt{n!}}\|u\EXP(f)\|
\sum_{\bar{k}(n)}
\sum_{\bar{\varepsilon}(n)}\|\delta(\bar{k}(n),\bar{\varepsilon}(n))
(x)\|, \]
since by  lemma \ref{lema1}(iii), the sum grows as $n$-th power,
 $ q_t(x)=0, \forall x\in
{\mathcal A}_{loc}, $
showing the uniqueness of the solution.
As $1 \in \mathcal A_{loc}$(in fact it is of empty support)
by QSDE \ (\ref{qsde}) it follows that $j_t(1)=1.$
\\(b)For $u\EXP(f),v\EXP(g)\in h\otimes\mathcal E(\mathcal C)$ and  $ x,y\in
\mathcal A_{loc}
 ,$  \    by induction we have, \[\langle
j_t^{(n)}(x^{*}) u\EXP(f),v\EXP(g) \rangle =\langle
u\EXP(f),j_t^{(n)}(x)v\EXP(g) \rangle\] Now as $ n $ tends to $ \infty $
we get \[ \langle j_t(x^{*}) u\EXP(f),v\EXP(g) \rangle =\langle
u\EXP(f),j_t(x)v\EXP(g) \rangle.\] Define
\[\Phi_t(x,y)= \langle u\EXP(f),j_t(xy)v\EXP(g) \rangle-\langle
j_t(x^{*}) u\EXP(f),j_t(y)v\EXP(g) \rangle  \]

Now for $l=1,2,\cdots,7 $~setting\\
$({{\zeta}_{k}}(l),{{\eta}_{k}}(l))=
(\delta_k,id),(id,\delta_k),(\delta_k^\dag,1),
(id,\delta_k^\dag),({\mathcal L}_{k},id),(id,{\mathcal L}_{k}) $ and
$ (\delta_k^\dag,\delta_k)$ \\ respectively,\ one has

\[|\Phi_t(x,y)| \leq c_{f,g}^n \sum_{l_n,\cdots,l_1} \int_0^t
\int_0^{s_{n-1}}\cdots\int_0^{s_1}\]
\be \label{homo2} \sum_{k_n,\cdots,k_1}|\Phi_{s_1}({{\zeta}_{k_n}}(l_n)\cdots
{{\zeta}_{k_1}}(l_1)x,\eta_{k_n}(l_n)\cdots \eta_{k_1}(l_1)y)|
ds_{0}\cdots ds_{n-1}, \ \forall n \geq 1,\ee
where $
c_{f,g}=({1+t_0}^{1/2})({\|f\|}_\infty+{\|g\|}_\infty).$
By quantum Ito formula and cocyle  properties of structure operators, i.e.
 $ {\hat{\mathcal L} }(xy)=x{\hat{\mathcal L}}(y) +
{\hat{\mathcal L} }(x)y+\sum _{k\in\MBBZ^{d}}\delta_k^\dag(x) \delta_k(y), $
we
have,
\bean \lefteqn{\Phi_t(x,y)}\\
&&={\int_0}^t\sum_{k} \{\Phi_s(\delta_k(x),y)
+\Phi_s(x,\delta_k(y))\}f_{k}(s)ds\\
&&+{\int_0}^t\sum_{k}
\{\Phi_s(\delta_k^\dag(x),y)+\Phi_s
(x,\delta_k^\dag(y))\}\bar{v}_k(s) ds\\
&&+{\int_0}^t\sum_{k} \{\Phi_s({\mathcal
L}_{k}(x),y)+\Phi_s(x,{\mathcal L}_{k}(y)) +
\Phi_s(\delta_k^\dag(x),\delta_k(y))\}ds \eean

which gives the estimate for $ n=1,$
\be \label{homo3} |\Phi_t(x,y) |\leq c_{f,g} \sum_{l=1\cdots
7}{\int_0}^t\sum_{k}|\Phi_s({{\zeta}_{k}}(l)(x),{{\eta}_{k}}(l)(y))|ds \ .
\ee  Now if we   assume \
(\ref{homo2}) for some $ m>1,$ application of \
(\ref{homo3}) gives  the required estimate  for $ n=m+1. $

Before going to further estimate of $  |\Phi_t(x,y) |, $  by \
(\ref{exist3}) ,  \ (\ref{exist4}) and lemma \ \ref{lema1} (iv),
note  the following,
\\(1)\ For any  $ n$-tuple $ (l_1,l_2\cdots l_n)$ in $\{1,2 \cdots 7\}$ \
\[ \sum_{k_n,\dots k_1}\|j_{s}({{\zeta}_{k_n}}(l_n)\cdots
{{\zeta}_{k_1}}(l_1)(x) \cdot\eta_{k_n}(l_n)\cdots
\eta_{k_1}(l_1)(y)) v\EXP(g)\|\] \be \label{homo21}\leq C_{g,x,y} \{
(1+\|r\|) (1+2\vartheta_1(r) c_{x,y})\}^{2n} \|v\EXP(g)\| \ee where for any
$ g\in \mathcal C $

\[C_{g,x,y}=1+\sum_{m\geq 1}3^m \frac{(t_0c_{g})^{m/2}}{\sqrt{m!}}{\{(1+\|r\|)
(1+2\vartheta_1(r) c_{x,y})\}}^{2m}.\]
(2)\ For any
 $ s\leq t_0,\ p\leq n $ and $ \bar{\varepsilon}(p) ,$  \

\[\sum_{\bar{k}(p)}\|j_{s}\{
\delta(\bar{k}(p),\bar{\varepsilon}(p))(y)\}v\EXP(g)\|\] \be
\label{homo22} \leq C_{g,x,y} {\{(1+\|r\|) (1+2\vartheta_1(r)
c_{x,y})\}}^n \|v\EXP(g)\| \ee
(3)\  Since $\vartheta_p(x)=\vartheta_p(x^{*})$ and $
\{\delta(\bar{k}(p),\bar{\varepsilon}(p))(x)\}^{*} $ can also be
written as $ \delta(\bar{k}(p),\bar{\varepsilon}^\prime(p))
(x^{*}) $ for some $ \bar{\varepsilon}^\prime(p) ,$  \   we have
\[ \sum_{\bar{k}(p)}\|j_{s}
\{\delta(\bar{k}(p),\bar{\varepsilon}(p))(x)\}^{*}u\EXP(f)\|\]
\be \label{homo23} \leq  C_{f,x,y}{\{(1+\|r\|) (1+2\vartheta_1(r)
c_{x,y})\}}^n \|u\EXP(f)\|
 \ee
 Now for any fixed $n$-tuple $(l_1,\cdots,l_n)  $ consider,
\[\sum_{\bar{k}(n)}|\Phi_s({{\zeta}_{k_n}}(l_n)\cdots
{{\zeta}_{k_1}}(l_1)x,\eta_{k_n}(l_n)\cdots \eta_{k_1}(l_1)y)|,\]
by definition of $ \Phi_s,$ it is

\[ \le  \sum_{k_n,\dots k_1}\|u\EXP(f)\|\cdot \| j_{s}({{\zeta}_{k_n}}(l_n)\cdots
{{\zeta}_{k_1}}(l_1)x \cdot\eta_{k_n}(l_n)\cdots \eta_{k_1}(l_1)y)
v\EXP(g)\|   \]
\[ + \| j_{s} \{ {({{\zeta}_{k_n}}(l_n)\cdots
{{\zeta}_{k_1}}(l_1)(x))}^{*}\} u\EXP(f)\|
+j_{s}(\eta_{k_n}(l_n)\cdots \eta_{k_1}(l_1)(y)) v\EXP(g)\|,\] now
the estimates \ (\ref{homo21}),\ (\ref{homo22}) and \
(\ref{homo23}) gives,
\bean \lefteqn{\sum_{\bar{k}(n)}|\Phi_s({{\zeta}_{k_n}}(l_n)\cdots
{{\zeta}_{k_1}}(l_1)x,\eta_{k_n}(l_n)\cdots \eta_{k_1}(l_1)y)|}\\
&&
\leq {\{ (1+\|r\|) (1+2\vartheta_1(r) c_{x,y})\} }^{2n} \|u\EXP(f)\|\cdot
\|v\EXP(g)\| ( C_{g,x,y}+ C_{f,x,y} C_{g,x,y} )\\
&&= C {\{(1+\|r\|) (1+2\vartheta_1(r)
c_{x,y})\} }^{2n} \eean with $ C=\|u\EXP(f)\|\cdot \|v\EXP(g)\| (
C_{g,x,y}+ C_{f,x,y} C_{g,x,y}),$ now by \ (\ref{homo2}),
\[|\Phi_t(x,y)|\leq  C \ \frac{(7 \ t_0c_{f,g})^n}{n!} \{ (1+\|r\|) (1+2\vartheta_1(r)
c_{x,y})\}^{2n},~\forall~n\geq 1.\] so $ \Phi_t(x,y)=0.$
\\(c) Let $\xi= \sum c_j u_j \EXP(f_j) $ (vector in algebraic tensor product  of $
\MBFH_0$ and $  \mathcal E(\mathcal C)).$ If $y
 \in  {\mathcal A}_{loc}^+,$  $ y $ is actually an
$ N^{|y|}\times N^{|y|} $-dim positive matrix and hence it admits
a unique  square root $ \sqrt{y} \in {\mathcal A}_{loc}^+.
 $ For any $ x\in {\mathcal
A}_{loc}^+, $ setting $ y=\sqrt{\|x\|1-x} $ so that $ y\in {\mathcal
A}_{loc}^+, $  we get
\bean \lefteqn{\|j_t(y)\xi\|^2=\langle j_t(y)\xi,j_t(y)\xi\rangle}\\
&&=\sum  \bar{c_i}c_j \langle j_t(y)u_i  \EXP(f_i),j_t(y)u_j \EXP(f_j)\rangle\\
&&=\sum  \bar{c_i}c_j \langle u_i \EXP(f_i),j_t(\|x\|1-x)u_j \EXP(f_j)\rangle \
\mbox{(by (b))}\\
&&= \| x \|\cdot {\| \xi \|}^2 - \langle \xi,j_t(x)\xi\rangle\eean
( where we have used the fact that $1\in {\mathcal A}_{loc} $ and
 $ j_t(1)=1 $ ) . Now let $ x\in {\mathcal
A}_{loc} $ be arbitrary and applying the above for $ x^*x $ and
by (b) we get,
\bean \lefteqn{\|j_t(x)\xi\|^2=\langle j_t(x)\xi,j_t(x)\xi\rangle}\\
&&=\sum  \bar{c_i}c_j \langle j_t(x)u_i \EXP(f_i),j_t(x)u_j \EXP(f_j)\rangle\\
&&=\sum  \bar{c_i}c_j \langle u_i \EXP(f_i),j_t(x^*x)u_j  \EXP(f_j)\rangle \\
&&= \langle \xi,j_t(x^*x)\xi\rangle\\
&&\le  \| x^*x \|\cdot {\| \xi \|}^2= \| x \|^2 \cdot {\| \xi \|}^2\\
&&\mbox{or}\ \|j_t(x) \xi \| \le \| x \| \cdot \| \xi \|.\eean This inequality
obviously extends to all $\xi\in   \MBFH_0 \otimes \Gamma. $ Noting
 that
$j_t(1)=1,\ \forall t,$  we get \[ \|j_t(x)\|\le \| x \| \
\mbox{and} \ \|j_t\|= 1. \] Thus $j_t$ extends uniquely to  a
unital $C^*$-homomorphism satisfying QSDE \ (\ref{qsde}) and
hence is a EH flow on $ \mathcal A $ with $ P_t $ as its
expectation semigroup.

\end{proof}

We have also obtained  an EH type dilation for the CP semigroup $ P_t^\phi$
associated with the partial state  $\phi_0$. Note that the generator
$\hat{\mathcal L^\phi}$ of $ P_t^\phi$ satisfies \[\hat{\mathcal L^\phi}
(x)= \sum_{k\in
 \MBBZ^d} \frac{1}{2}\sum_{m=1}^{N^\prime} [{ L_k^{(m)}}^*,x]L_k^{(m)}+
 { L_k^{(m)}}^*[x,L_k^{(m)}],\forall x\in \A_{loc}.\]
  Now we have the following,

\bthm
Let $\hat{\mathcal L^\phi}$ and $ P_t^\phi$ as discussed earlier, then
\\(a) For each  $k\in
 \MBBZ^d $ and $t\ge 0 $ there exist  unique solution $\eta_t^{(k)} (x_{(k)})$
 for the  QSDE,
 \be \label{eh-ergodic}
d\eta_t^{(k)}(x)=\eta_t^{(k)}
(\sum_{m=1}^{N^\prime}[{L_k^{(m)}}^*,x_{(k)}])
d{a_k}(t)+\eta_t^{(k)}
(\sum_{m=1}^{N^\prime}[x_{(k)},L_k^{(m)}])
da_k^\dag(t)
+\eta_t^{(k)}(\mathcal L_k^\phi x_{(k)})dt,\ee
\[{j_0}(x_{(k)})=x_{(k)} \otimes 1_{\Gamma},\  \forall x_{(k)} \in
\mathcal A_k  \] and $\eta_t^{(k)}$ is a  unital  $*$-homomorphism
 from $\mathcal A_k  $
in to $\mathcal A_k \otimes \mathcal B (\Gamma).$
Moreover, for different $k $ and $k^\prime, \eta_t^{(k)}$ and
 $\eta_t^{({k^\prime})}$  commute in the sense that,
$\eta_t^{(k)}(x_{(k)})$ and $\eta_t^{({k^\prime})}(x_{k^\prime})$  commute for every
 $ x_{(k)} \in
\mathcal A_k $ and $x_{k^\prime} \in
\mathcal A_{k^\prime},$
\\(b)There exist unique  unital  $*$-homomorphism
$\eta_t$  from $\mathcal A_{loc} $ in to $\mathcal A \otimes \mathcal B (\Gamma)$
such that it coincide  with $\eta_t^{(k)}$ on $\mathcal A_k, $
\\(c)The $\eta_t$ extends uniquely as a unital  $C^*$-homomorphism
 from $\mathcal A $ in to $\mathcal A \otimes \mathcal B (\Gamma).$
\ethm

\begin{proof}
(a)
 For any  $k\in
 \MBBZ^d, t\ge 0 $ and $ x_{(k)} \in
\mathcal A_k, $  consider the  QSDE\
 (\ref{eh-ergodic}). Here  we have only finitely many nontrivial structure maps
 on the unital $C^*$-algebra $\mathcal A_k,$
 satisfying structure equation . So
 there exist a unique solution $\eta_t^{(k)}(x_{(k)})$ and
  $\eta_t^{(k)}$ is a unital  $*$-homomorphism
 from $\mathcal A_k  $
in to $\mathcal A_k \otimes \mathcal B (\Gamma).$
Note that for  different $k $ and $ k^\prime$  associated
  structure maps are commuting. Hence,
$\eta_t^{(k)}(x_{(k)})$ and $\eta_t^{({k^\prime})}(x_{(k^\prime)})$  commute for every
 $ x_{(k)} \in
\mathcal A_k $ and $x_{(k^\prime)} \in
\mathcal A_{k^\prime}.$\\

\noindent (b)Now for any finite $\Lambda\subseteq \MBBZ^d,t\ge0 $
 and simple tensor element
$x_\Lambda=\prod_{k \in \Lambda}x_{(k)} \in
\mathcal A_\Lambda, $  if we set
\[\eta_t^{(\Lambda)}(x_\Lambda)= \prod_{k \in \Lambda}\eta_t^{(k)}(x_{(k)})\]then
$\eta_t^{(\Lambda)}$ is a well
 defined map on $\mathcal A_\Lambda$ to $\mathcal A_\Lambda\otimes \mathcal B
 (\Gamma)$ due to the fact that $\eta_t^{(k)}$'s
  commute. Differentiating
$\eta_t^{(\Lambda)}(x_\Lambda)$ with respect to $t,$ it follows that
$ \eta_t^{(\Lambda)}(x_\Lambda)$ satisfies the QSDE,
\be \label{ehergodic}
d\eta_t^{(\Lambda)}(x_\Lambda)=\sum_{k\in \Lambda}\eta_t^{(\Lambda)}
(\sum_{m=1}^{N^\prime}[{L_k^{(m)}}^*,x_\Lambda])
d{a_k}(t)+\sum_{k\in \Lambda}\eta_t^{(\Lambda)}
(\sum_{m=1}^{N^\prime}[x_\Lambda,L_k^{(m)}])
da_k^\dag(t)
+\eta_t^{(\Lambda)}(\mathcal L_k^\phi x_\Lambda)dt,\ee
\[\eta_0^{(\Lambda)}(x_\Lambda)=x_\Lambda \otimes 1_{\Gamma}. \]
In order to show,  \be \label{simlehomo} \eta_t^{(\Lambda)}(xy)=
\eta_t^{(\Lambda)}(x)\cdot \eta_t^{(\Lambda)}(y),\ \mbox{for every
simple tensor elements }\  x,y\in \mathcal A_{loc},\ee
without loss of generality (since each $\eta_t^{(k)}$'s  are unital, for
 finite subsets $\Lambda \subseteq \Lambda^\prime ,\eta_t^{(\Lambda^\prime)}$ agree
 with
$\eta_t^{(\Lambda)}$ for   simple tensor elements in $\mathcal A_\Lambda$) assume
$ x,y\in \mathcal A_\Lambda $ for some
 finite  $\Lambda\subseteq \MBBZ^d$ so that
$x=\prod_{k \in \Lambda}x_{(k)} \in
\mathcal A_\Lambda $ and $y=\prod_{k \in \Lambda}y_{(k)} \in
\mathcal A_\Lambda$ with identity  component out side their respective  support.
 Now consider, \bean \lefteqn{\eta_t^{(\Lambda)}(xy)=
\eta_t^{(\Lambda)}\prod_{k \in \Lambda}(x_{(k)} y_{(k)})
=\prod_{k \in \Lambda}\eta_t^{(k)}(x_{(k)} y_{(k)})}\\
&&=\prod_{k \in \Lambda}\eta_t^{(k)}(x_{(k)})\eta_t^{(k)}( y_{(k)})=
\prod_{k \in \Lambda}\eta_t^{(k)}(x_{(k)})\prod_{k \in \Lambda}\eta_t^{(k)}( y_{(k)}).\eean
Thus,\ (\ref{simlehomo}) follows. Similarly for $x=\prod_{k \in \Lambda}x_{(k)},$
\be \label{star}  \eta_t^{(\Lambda)}(x^*)=(\eta_t^{(\Lambda)}(x))^*.\ee
Now we define $ \eta_t $ on $\mathcal A_{loc}$ as follows,
note that any element $x\in \mathcal A_{loc} $ can be written as a
linear combination of
simple tensor elements $\{U_g:g\in \mathcal G\}, x=\sum_{g\in \G }c_gU_g$
with $c_g=o$  when $supp(g)$ is  outside the $supp(x)=\Lambda,$   set
\[\eta_t(x)=\sum_{g\in \G }c_g \eta_t^{(\Lambda)}(U_g)\]
Let $  x $ and $y\in \mathcal A_{loc}, x=\sum_{g\in \G }c_gU_g$ and
 $ y=\sum_{h\in \G }c_hU_h$ such that $supp(x)=supp(y)=\Lambda,$  consider
\bean \lefteqn{\eta_t(xy)
=\eta_t(\sum_{g,h\in \G }c_g c_hU_gU_h)}\\
&&=\sum_{g,h\in \G }c_g c_h\eta_t^{(\Lambda)}(U_gU_h)=
\sum_{g,h\in \G }c_g c_h\eta_t^{(\Lambda)}(U_g)\eta_t^{(\Lambda)}
(U_h)\ \mbox{(by (\ref{simlehomo})  )}\\
&&=\eta_t(\sum_{g\in \G }c_gU_g)\eta_t(\sum_{h\in \G }c_hU_h).\eean
So $\eta_t(xy)=\eta_t(x)\eta_t(y)$ and by \ \ref{star} it follows that
$\eta_t(x^*)=(\eta_t(x))^*,\forall x\in \mathcal A_{loc}.$ Thus
$\eta_t$ is a  unital  $*$-homomorphism  from $\mathcal A_{loc}  $
in to $\mathcal A \otimes \mathcal B (\Gamma).$
(c)\ (Proof is same  as that of theorem  \ \ref{theorem1}(c))\\
Let $ x\in {\mathcal
A}_{loc} $ then  $ \| x \|^21-x^*x \in {\mathcal
A}_{loc}^+$ (in fact it is belong to some finite dimensional matrix algebra
$\mathcal
A_\Lambda $)   so $\sqrt{\| x \|^21-x^*x} \in \mathcal
A_{loc}^+.$ Since $\eta_t$ is a  unital  $*$-homomorphism on
$\mathcal A_{loc},$
\[\eta_t(\| x \|^21-x^*x )\ge 0\]
\[\Rightarrow \eta_t(x^*x )\le \| x \|^21\]
\[\Rightarrow \|\eta_t(x^*x )\|\le \| x \|^2\]
\[\Rightarrow \|\eta_t(x)\|\le \| x \|,\]
So $\eta_t$ extends uniquely as a
unital $C^*$-homomorphism  from $\mathcal A  $
in to $\mathcal A \otimes \mathcal B (\Gamma).$


\end{proof}

\section{Covariance  of EH flow }
In this section  we shall prove that the Evans-Hudson flows
constructed in the last section  is covariant. Let  $\mathcal B $ be a $
C^* $ ( or von Neumann) algebra, $ G $ be a locally compact group equipped
with an action $ \alpha $ on $ \mathcal B .$\  Let $ \{T_t:t>0\} $ be a
covariant  CP semigroup on $ \mathcal B $ w.r.t.  $\alpha,$ that is,
$$\alpha_{g}\circ T_t (x)=T_t\circ \alpha_{g}(x),\forall t\ge
0,g\in G , x\in \mathcal B .$$ Then a natural question arises ,\ does
there exist  a covariant EH dilation for $ \{T_t\}. $ The question is discussed
in \cite{CGS} for uniformly continuous CP semigroup.\ There is  no
such general result for  CP semigroups with unbounded generators.\

We shall show  the EH flow $ \{j_t\}  $ and $ \{\eta_t\}$ constructed in the
previous section is covariant w.r.t. the actions $ \tau $ and $\lambda $( $\lambda
$  to be introduced later in this section)
of the  group $ \MBBZ^{d} $ .\\
It can be easily observed that \be \label{cov0}\delta_k
\tau_j=\tau_j\delta_{k-j} \mbox{ and}\
{\delta^\dag}_k\tau_j=\tau_j{\delta^\dag}_{k-j}, \ \forall j,k\in
\MBBZ^d \ee and we have the following lemma,
 \blema \label{cov1} (i)$ \hat{\mathcal
L}\tau_j(x)=\tau_j\hat{\mathcal L}(x),\
\forall x\in Dom( \hat{\mathcal L}),$\\
(ii)$ P_t\tau_j=\tau_j P_t, $ i.e. $ P_t $ is covariant.
 \elema
 \begin{proof}
(i) Note that $\C^1(\A)$ is $\tau $ invariant and thus for
 $ x\in {\mathcal C}^1({\mathcal A}),$
 \bean \lefteqn{ {\mathcal
L}(\tau_j(x))= \frac{1}{2}\sum_{k \in
\MBBZ^d}\delta_k^\dag(\tau_j(x))r_k+
r_k^*\delta_k(\tau_j(x))}\\
&&=\frac{1}{2}\sum_{k \in
\MBBZ^d}\tau_j\delta_{k-j}^\dag(x) r_k+
{r_k}^{*}\tau_j\delta_{k-j}(x) ~~~~ \mbox{(by \ref{cov1} )}\\
&&=\frac{1}{2} \tau_j \{\sum_{k \in
\MBBZ^d}\delta_{k-j}^\dag(x) r_{k-j}+
r_{k-j}^{*}\delta_{k-j}(x)\}\\
&&=\tau_j(\mathcal L(x)).\eean
For $ x\in  Dom(\hat{\mathcal
L}),$ choose a sequence $\{x_n\} $ in $ {\mathcal
C}^1({\mathcal A})$ and $ y\in \mathcal A $ such that
$y=\hat{\mathcal
L}(x) $  and $ x_n $ and ${\mathcal L}(x_n) $ converges to $ x $  and
 $ y $ respectively. Now, for any $ j\in \MBBZ^d$ applying the automorphism
 $ \tau_j, \tau_j(x_n) $  and
$ \tau_j{\mathcal
L}(x_n)$ converges to $ \tau_j(x)$ and
 $\tau_j(y)$ respectively.
Since $ x_n\in {\mathcal C}^1({\mathcal A}),
 {\mathcal L}(\tau_j(x_n))=\tau_j{\mathcal L}(x_n) $  and we get
\[\tau_j(x)\in Dom(\hat{\mathcal L}) \ \mbox{and }\ \hat{\mathcal
L}\tau_j(x)=\tau_j\hat{\mathcal L}(x).
\]
(ii)By (i), for $ x\in Dom(\hat{\mathcal L})$\ and $\  0\leq s\leq t$ we have,
\[\frac{d}{ds}P_s \circ \tau_j \circ P_{t-s}(x)=P_s \circ \hat{\mathcal L}\circ \tau_j
\circ P_{t-s}(x)-P_s \circ\tau_j \circ \hat{\mathcal L} \circ
P_{t-s}(x)=0\]
This implies that $ P_s \circ \tau_j \circ P_{t-s}(x) $ is independent
of $ s $ for every $ j $ and $ 0\leq s\leq t.$  Setting $ s=0 $ and $ t $
respectively  and using the fact that $ P_t $ is bounded  we get
$   P_t\tau_j=\tau_j P_t.$
\end{proof}

Note that $j_t :\mathcal A \rightarrow
{\mathcal A} \otimes {\mathcal B}(\Gamma(L^2(\MBBR_+,\MBFK_0))),$
where $\MBFK_0 =l^2(\MBBZ^d)$ with canonical basis
$\{e_k\},$ as mentioned earlier. Define  the canonical bilateral shift
 $ s $ by  $ s_{j}e_{k}=e_{k+j},\forall j,k\in \MBBZ^{d} $ and let
$\gamma_{j}=\Gamma( 1\otimes s_{j}),$ the second quantization of
$ 1\otimes s_{j} $ i.e. $ \gamma_j \EXP(\sum f_l(.)e_l)=\EXP(\sum
f_l(.)e_{l+j}), $ this defines a unitary representation of
$\MBBZ^{d} $ in $ \Gamma $ and further we set $\sigma={\tau}
\otimes \lambda  $  on $ {\mathcal A} \otimes {\mathcal B}(\Gamma)
$ where $ \lambda_j(y)=\gamma_j  y \gamma_{-j}, \ \forall y\in
{\mathcal B}(\Gamma). $

By definition of fundamental processes $ a_k(t): a_k(t)\EXP(g)=\int_0^t
g_k(s)ds \ \EXP(g),$ observed that
 \bean \lefteqn{\lambda_j a_k(t)\EXP(g)=\gamma_j
a_k(t) \gamma_{-j}\EXP(g)=\gamma_j a_k(t)\EXP(\sum \langle
g,e_{l+j}\rangle (\cdot) e_l)}\\
&&= \int_0^t\langle g,e_{k+j}\rangle (s)ds \ \gamma_j
(\EXP(\sum \langle g,e_{l+j}\rangle (\cdot)e_l)\\
&&= \int_0^t\langle g,e_{k+j}\rangle (s)ds \
(\EXP(\sum \langle g,e_{l+j}\rangle (\cdot)
e_{l+j})\\
&&=a_{k+j}(t)\EXP(g)\eean
and since \[ \langle \EXP(f),\lambda_j a_k(t)\EXP(g)\rangle=\langle\lambda_j
 a_k^\dag(t)\EXP(f),\EXP(g)\rangle,\]
it follows that
\be \label{aa+} \lambda_j a_k(t)=a_{k+j}(t) \ \mbox{and}\
\lambda_j a_k^\dag (t)=a_{k+j}^\dag (t).\ee

\bthm The Evans-Hudson flow $j_t$ of the CP semigroup $P_t$ is
covariant,i.e. \[\sigma_j j_t\tau_{-j}(x)=j_t(x),\forall x\in
{\mathcal A},\ t \ge 0,k\in \MBBZ^d. \] \ethm
\begin{proof}
For  a fixed $ j
\in \MBBZ^d,$ set  $j_t^\prime=\sigma_j j_t\tau_{-j},\ \forall t\ge 0. $
Using QSDE \ (\ref{qsde})\ and lemma \ \ref{cov1},\ (\ref{cov0}),\ (\ref{aa+})
we get for $ x\in \mathcal A_{loc}, $

\bean \lefteqn{j_t^\prime(x)-j_0^\prime (x)}\\
&&=\int_0^t\{ \sum_{k\in \MBBZ^d}\sigma_j j_s
(\delta _k^\dag (\tau_{-j}(x)))
da_k(s)+\int_0^t\sum_{k\in \MBBZ^d}\sigma_j j_s(\delta _k (\tau_{-j}(x)))
da_k^\dag(s)\\
&&~~~~~~~+\int_0^t \sigma_j j_s({\hat{\mathcal L}} (\tau_{-j}(x)))ds\\
&&=\int_0^t \sum_{k\in \MBBZ^d}\sigma_j j_s\tau_{-j}
(\delta_{k+j}^\dag (x))
d{a_{k+j}}(s)+\int_0^t \sum_{k\in \MBBZ^d}\sigma_j j_s\tau_{-j}(\delta_{k+j}(x))
d{a_{k+j}^\dag}(s)\\
&&~~~~~~~+\sigma_j j_s\tau_{-j}({\hat{\mathcal L}} (x))ds \\
&&=\int_0^t\{\sum_{k\in \MBBZ^d}j_s^\prime(\delta _k^\dag (x))
d{a_k}(s)+\int_0^t\sum_{k\in \MBBZ^d} j_s^\prime(\delta_k (x))
d{a_k^\dag}(s)+\int_0^t j_s^\prime({\hat{\mathcal L}} x)ds.\\
&& \mbox{Since,}\  j_0^\prime(x)=\sigma_j j_0\tau_{-j}(x)=
\sigma_j (\tau_{-j}(x)\otimes 1_\Gamma) =x\otimes 1_\Gamma=j_0(x),\eean
$j_t^\prime(x)=j_t(x)$ for all $ t\ge 0 $ and $ x\in \A_{loc},$ by uniqueness  of QSDE \
(\ref{qsde}).  As  both $ j_t^\prime $
and $ j_t $ are bounded maps, it follows that $j_t^\prime = j_t.$
\end{proof}

\brmrk By similar argument as  above, the EH flow   for the CP semigroup
 $P_t^\phi,$ can be seen to be
 covariant with respect to the same actions.
 \ermrk
\section{Ergodicity of the EH flows}

Recall the QDS  $P_t^\phi$ associated with the  partial state $\phi_0,$ for which
we have constructed
EH flows $\eta_t$ in section 3. Note that
 $P_t^\phi$  has a unique ergodic state $\Phi$.
We have the following result on ergodicity  of $\eta_t$
w.r.t. the weak topology.
\bthm
 The EH flow $\eta_t$ of CP semigroup $P_t^\phi$
has also the  unique ergodic state  $\Phi$, in a sense that
\[\eta_t(x)\rightarrow\Phi(x)1_\Gamma \ \mbox{weakly}\
\forall x\in {\mathcal A}.\]
\ethm
\begin{proof} Since $\eta_t$ and $P_t^\phi$ are norm contractive, $\A_{loc}$
is norm-dense in $\A,$    and $ P^\Phi_t(x)$ converges to
$\Phi(x)\otimes 1$ for all $x\in \A,$ it is enough to show that $\eta_t(x)-P^\Phi_t(x) \otimes 1
 \rightarrow 0$ weakly  as $t \rightarrow \infty.$ Furthermore, it suffices to show that
  $ \langle \xi_1, (\eta_t(x) -P_t^\Phi(x) \otimes 1) \xi_2 \rangle \rightarrow 0$ as $t\rightarrow \infty$,
   where $\xi_1,\xi_2$  vary over the linear span of vectors of the form $ve(f)$, with
    $f=\sum_{|k|\leq n}f_k\otimes e_k$ for some $n$ and $f_1,...f_n \in L^1(\MBBR_+)
    \bigcap L^2(\MBBR_+).$

For notational simplicity denoting  the bounded derivations on $\A,$
\[x \mapsto \sum_{m=1}^{N^\prime}[x,L_k^{(m)}]
\ \mbox{and}\ x \mapsto \sum_{m=1}^{N^\prime} [{ L_k^{(m)}}^*,x] \]
by $\rho_k$  and $\rho_k^\dag$ respectively,
note that $\eta_t$ satisfies the QSDE,
\be \label{ehergodic1}
d\eta_t(x)=\sum_{k\in \MBBZ^d}\eta_t
(\rho_k^\dag(x))
d{a_k}(t)+\sum_{k\in \MBBZ^d }\eta_t
(\rho_k(x))da_k^\dag (t)
+\sum_{k\in \MBBZ^d }\eta_t(\mathcal L_k^\phi ( x))dt,\ee
\[{\eta_0}(x)=x \otimes 1_{\Gamma},\forall x \in \A_{loc}. \]
For  $t\ge 0,u,v\in\MBFH_0 $ and $ f,g\in
L^2(\MBBR_+,\MBFK_0)\bigcap L^1(\MBBR_+,K_0) $ such that $
f=\sum_{|k|\leq n}f_k\otimes e_k$ \ and $g=\sum_{|k|\leq
n}g_k\otimes e_k$ and $x\in \mathcal A_{loc},$
consider the following,
\bean \lefteqn{|\langle u\EXP(f),[\eta_t(x)-P_t^\phi(x)\otimes 1]v\EXP(g)\rangle|}\\
&& = |\langle u\EXP(f),[\int_0^t\sum_{k \in \MBBZ^d}\eta_q\{ \rho_k
(P_{t-q}^\phi(x))\}da_k^\dag (q) + \eta_q\{\rho_k^\dag
(P_{t-q}^\phi(x))\}da_k(q)]v\EXP(g)\rangle| \\
&& \leq  \sum_{|k|\leq n}\int_0^t | \langle u\EXP(f),\eta_q\{\rho_k
(P_{t-q}^\phi(x))\}v\EXP(g)\rangle |\  \|g(q)\|dq \\
&&+\sum_{|k|\leq n}\int_0^t | \langle u\EXP(f),\eta_q\{\rho^\dag _k
(P_{t-q}^\phi(x))\}v\EXP(g)\rangle |\ \|f(q)\|dq \eean
 As $\eta_t,\  P_t^\phi$ are  contractive and
$P_t^\phi(x)$ tend to $\Phi(x)1$ as $t$ tend to $\infty,$ further
$\rho_k
$ and $\rho_k^\dag $ are   uniformly  bounded and $\rho_k(1)=\rho_k^\dag(1)=0$
for all $k\in \MBBZ^d,$ we have,

$| \langle u\EXP(f),\eta_q\{\rho_k
(P_{t-q}^\phi(x))\}v\EXP(g)\rangle |$ and $ | \langle u\EXP(f),\eta_q\{\rho^\dag _k
(P_{t-q}^\phi(x))\}v\EXP(g)\rangle |
\leq M, $ \\ for
 some constant $ M $ independent of $t $ and $ q.$
So since $ f,g\in L^1(\MBBR_+,K_0),$  both the terms
of the above expression tend to $0$
as $t $ tends to $\infty.$

\end{proof}
\brmrk The $\eta_t(x)$ \ does not converge  strongly, for if did, then $x\mapsto\Phi(x)
\otimes 1_\Gamma $ would be a homomorphism, i.e. $\Phi$ would be a
multiplicative non zero functional on the UHF algebra $\A,$ contradictory
to the fact that $\A$ does not have any such functional. \ermrk
\brmrk
If we look at the  perturbation of the QDS  $P_t^\phi$
by QDS associated with single supported  $r\in \A_0,$ then
by the same  argument   used  in  the construction of EH dilation  for
the unperturbed semigroup will go through  and one
can obtain an EH dilation  for the  perturbed one. For small perturbation
 parameter $c\ge 0$
for which  a unique ergodic state  exists, the  EH flow also
admits the   same  unique ergodic state in the  above sense.
 \ermrk


\begin{thebibliography}{CGS}
\bibitem{CGS}Chakraborty, P. S., Goswami, D. and  Sinha, K.  B.;
A covariant quantum stochastic dilation theory. Stochastics in finite and infinite
dimensions, 89--99, Trends Math., Birkhäuser Boston, Boston, MA, 2001.

\bibitem{GS1} Goswami, D. and Sinha, K.\ B.; \  Hilbert modules and stochastic
 dilation of a quantum dynamical semigroup on a von Neumann algebra, Comm. Math.
  Phys.  (1999)\ {\bf 205} no. 2, 377--403.

\bibitem{GS2} Goswami, D. and Sinha, K.\ B.; Stochastic dilation of symmetric
completely  positive   semigroups , submitted to J. of Funct. Anal.
\bibitem {GPS} Goswami, D., Pal, A.\ K.\ and Sinha, K.\ B.; {\it Stochastic
dilation of a quantum dynamical semigroup on a separable unital $C\sp *$-algebra}
. Infin. Dimens. Anal. Quantum Probab. Relat. Top. (2000)\ {\bf 3}, no. 1, 177--184.
\bibitem {HP} Hudson, R. L. and Parthasarathy, K. R.;
 Quantum Ito's formula and stochastic evolutions.
 Comm. Math. Phys.  (1984) {\bf 93}, no. 3, 301--323.
\bibitem {Ma} Matsui, Taku; Markov semigroups on UHF algebras. Rev. Math. Phys.
 (1993)\ {\bf 5,} no. 3, 587--600.
\bibitem{Mey} Meyer, P.A.;`` Quantum Probability  for Probabilist'' 2nd ed,
Lecture Notes in mathematics, Vol.1538, springer-Verlag, Heidelberg 1993.
\bibitem{Mo} Mohari, A.;  Quantum stochastic differential equations with
unbounded coefficients and dilations of Feller's minimal solution. Sankhy\=a Ser.
 A  (1991)\ {\bf 53,} no. 3, 255--287.
 \bibitem{MS}Mohari, A. and Sinha, K. B.; Stochastic dilation of minimal quantum
 dynamical
semigroup. Proc. Indian Acad. Sci. Math. Sci. (1992)\ {\bf 102,} no. 3, 159--173.

\bibitem {KRP} Parthasarathy, K.\ R.;\ `` An introduction to quantum stochastic
calculus'', Monographs in Mathematics, {\bf 85}, Birkh$\ddot{a}$user Verlag, Basel,
 1992.


\end{thebibliography}
\end{document}